\documentclass[letterpaper, 10 pt, conference]{ieeeconf}  
\IEEEoverridecommandlockouts
\usepackage{comment}
\usepackage{cite}
\usepackage{amsmath,amssymb,amsfonts}
\usepackage{graphicx}
\usepackage{textcomp}
\usepackage{color}
\usepackage{algorithmicx}
\usepackage{wrapfig}
\usepackage{float}

\usepackage[nolist]{acronym}
\usepackage{xcolor}
\usepackage{amsmath}
\usepackage{amssymb}
\usepackage{mathtools}
\usepackage{subcaption}
\usepackage{hyperref}
\hypersetup{colorlinks=true, linkcolor=black}
\usepackage[ruled]{algorithm2e}
\usepackage{mathrsfs}
\definecolor{ccolor}{RGB}{203,96,21}

\DeclarePairedDelimiterX{\inp}[2]{\langle}{\rangle}{#1, #2}

\DeclareMathOperator*{\argmin}{arg\!\,min}

\newtheorem{theorem}{Theorem}

\newtheorem{remark}{Remark}

\newtheorem{definition}{Definition}
\newcommand{\beq}{\begin{equation}}
\newcommand{\eeq}{\end{equation}}
\newcommand{\mat}[1]{{\bm{#1}}}

\newcommand{\greekvect}[1]{\boldsymbol{#1}}

\usepackage{bm}
\usepackage{bbm}

\begin{document}

\title{\LARGE \bf
Challenges in  Model Agnostic Controller Learning for   Unstable Systems}
\author{Mario Sznaier$^{1}$, Mustafa Bozdag$^{1}$
\thanks{This work was partially supported by NSF grant CNS--2038493, AFOSR grant FA9550-19-1-0005, ONR grant N00014-21-1-2431, and  the Sentry DHS Center of Excellence
under Award 22STESE00001-03-03}
\thanks{$^{1}$Robust Systems Lab,  ECE Department, Northeastern University, Boston, MA 02115. (bozdag.m@northeastern.edu, msznaier@coe.neu.edu)}}

\maketitle
\thispagestyle{empty}
\pagestyle{empty}

\begin{abstract} Model agnostic controller learning, for instance by direct policy optimization, has been the object of renewed attention lately, since it avoids a computationally expensive system identification step. Indeed, direct policy search has been empirically shown to lead to optimal  controllers in a number of cases of practical importance. However, to date, these empirical results have not been backed up with a comprehensive theoretical analysis for general problems.  In this paper we use a simple example to show that direct policy optimization is not directly generalizable to other seemingly simple problems. In such cases, direct optimization of a performance index can lead to unstable pole/zero cancellations, resulting in the loss of internal stability and unbounded outputs in response to arbitrarily small perturbations. We conclude the paper by analyzing several alternatives to avoid this phenomenon, suggesting some new directions in direct control policy optimization.

\end{abstract}

\section{Introduction}

Recently, there has been renewed interest in ``model free" control design techniques where the goal is to design a controller that optimizes performance based purely 
 on experimental data. These techniques are attractive since they hold the promise of optimizing performance while avoiding a computationally expensive systems identification step 
 \cite{sznaier2020}.  Indeed, direct control policy  optimization techniques have  achieved remarkable success in a range of classical control problems, ranging from Linear Quadratic (LQR, LQG) to $\mathcal{H}_\infty$ and controller auto-tuning.  While a general theory supporting these results is still emerging \cite{Hu2023}, recent results show that, in spite of lack of convexity, direct optimization can lead to optimal policies in these problems \cite{pmlr-v80-fazel18a,pmlr-v89-malik19a,fatkhullin2020optimizing,Kumar2021,Mehndiratta2021, Loquercio2022,Cheng2023,zheng2023benign} and provide bounds on the sample complexity \cite{Jovanovic2022,guo2023complexity}.  Hence, the hope is that these techniques can provide a viable alternative to the traditional systems identification-control design pipeline. 
The goal of this paper is to point out to the dangers of using direct optimization even in seemingly simple problems. As we illustrate with a simple first order system, direct policy optimization can lead to unstable pole zero cancellations and hence the loss of internal stability. In turn, this can result in unbounded signals in response to arbitrarily small perturbations to the control action.

The paper is organized as follows:
 Section \ref{sec:notation} introduces the notation and required definitions; Section \ref{sec:example} contains tha main result of the paper: a simple example where model agnostic performance optimization over all continuous stabilizing controllers leads to a pole/zero cancellation and loss of internal stability; Section \ref{sec:analysis}  connects this example with the empirical observation in \cite{Seiler2019} that adding noise during training increases robustness, and discusses some ideas to avoid the loss of internal stability. Section \ref{sec:conclusions} offers some conclusions and points out to directions for further research. 

\section{Notation and Definitions}\label{sec:notation}

$\|u\| \doteq \sqrt{\mat{u}^T\mat{u}}$ denotes the usual Euclidean norm in $R^n$. For a given
 sequence $x_k$, $\|x\|_2 \doteq \sqrt{\sum \|x_k\|^2}$. $\ell^2$ denotes the Hilbert space of real vector sequences with finite $\|x\|_2$, equipped with the inner product $\langle x,y \rangle \doteq \sum x_i y_i$. $\ell^\infty$ denotes the Banach space of bounded real sequences equipped with the norm $\|x_k\|_\infty \doteq \sup_k |x_k|$.
We will denote  with capital letters the  z-transform of sequences in $\ell^2$, e.g. 
$X(z) = \sum x_k z^{-k}$. By a slight abuse  of notation, sometimes we will write  $\|X(z)\|_2$ to denote the $\ell^2$ norm of the sequence $\{x_k\}$.  Given a sequence $x_k$,
we will denote by $(x_k)_\tau$ its truncation, that is,
$(x_k)_\tau = x_k$ for $0 \leq k \leq \tau$ and $(x_k)_\tau=0$ otherwise. In this context the extended space $\ell^\infty_e$ is defined as $\ell^\infty_e = \left \{u \colon (u)_\tau \in \ell^\infty \; \forall \; \tau \in [0, \infty)\right \}$. $\mathcal{RH}_\infty$ denotes the Lebesgue space of complex valued rational functions 
 with bounded analytic continuation
in $|z| > 1$, equipped with the norm 
$\|G(z)\|_{\mathcal{H}_\infty} \doteq \sup_{|\lambda|>1} |G(z)|$. In the sequel, we will represent a linear time-invariant system $\mathcal{G}:\ell^\infty_e\to \ell^\infty_e$ either by its convolution kernel $g$ or the transfer function $G(z)$. It is well known \cite{SpSzBook} that, if $\mathcal{G}$ is stable, then its  $\ell^2$ induced norm $\|\mathcal{G}\|_{\ell^2 \to \ell^2} = \|G(z)\|_{\mathcal{H}_\infty} $.
Finally, given a matrix $\mat{M}$, $\|\mat{M}\|_2$ denotes its $\ell^2 \to \ell^2$ induced norm.

\begin{definition}\cite{Khalil} An operator $H:\ell^\infty_e \to \ell^\infty_e$ is finite $\ell^\infty$-gain stable if there exist constants $\gamma \geq 0, \beta \geq 0$ such that
$\|(Hw)_\tau \|_\infty \leq \gamma \|(w)_\tau\|_\infty + \beta$,  $ \forall w \in \ell^\infty_e$ and  $\tau \in [0,\infty)$.
\end{definition}
In the sequel, by a slight abuse of notation we will restrict
this definition to the case where $\beta=0$, that is, we will only consider mappings where $H0=0$. Thus, in the case of linear systems, finite-$\ell^\infty$ gain stability reduces to the standard bounded-input bounded-output stability.
\begin{figure}[hbt]
\vskip -1em
\centering
\includegraphics[width=0.75\linewidth]{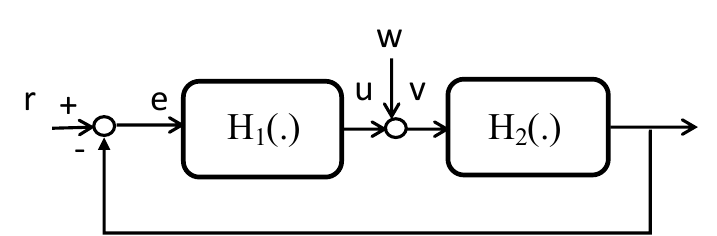} 
\caption{The closed loop is internally stable if all four mappings $\begin{bmatrix}r & w \end{bmatrix}^T \rightarrow \begin{bmatrix}e & v \end{bmatrix}^T$ 
are stable.}  \label{fig:internal} 
\end{figure}
\begin{definition} The loop shown in Fig. \ref{fig:internal} is finite $\ell^\infty$ gain   internally stable if all four mappings $\begin{bmatrix}r & w \end{bmatrix}^T \rightarrow \begin{bmatrix}e & v \end{bmatrix}^T$  are finite $\ell^\infty$ gain  stable.
\end{definition}

\section{A simple first order example}\label{sec:example}

\begin{figure}[H]
    \centering
     \vskip -0.5 cm
    \includegraphics[width=0.7\linewidth]{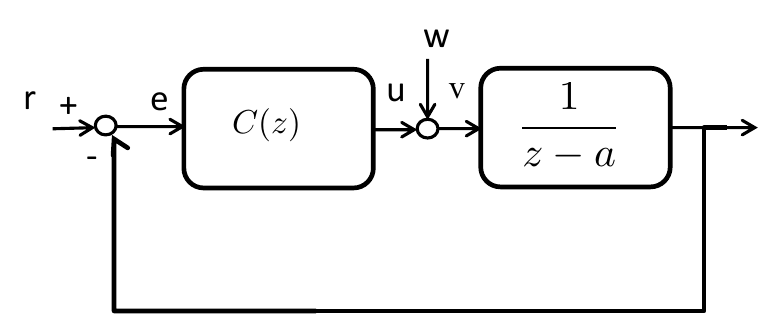} 
    \caption{Closed loop for the simple example.}  
    \label{fig:CL1} 
\end{figure}
\vskip 0.5 em

Here we present a simple example  where input/output optimization leads to loss of  internal stability: Given the open loop unstable system shown in Fig. \ref{fig:CL1}, where $a >1$, find a  controller $C$ that minimizes $\|e\|_2$ to an input $r$ of the form:
\beq \begin{small}
\label{eq:rectangular}
r_k = 
\begin{cases} 
r_o, & 0 \leq k \leq n-1 \\ 
0, & \text{otherwise} 
\end{cases}
\Rightarrow
R(z) = r_o \frac{z^n - 1}{(z-1)z^{n-1}}
\end{small}
\eeq
where both the amplitude $r_o$ and width $n$ are unknown and arbitrary. In the sequel, we will consider three scenarios: (A) optimization over all internally stabilizing LTI controllers; (B) optimization over all LTI controllers that only stabilize the mapping $r \to e$; and (C) optimization over all causal, time-invariant continuous nonlinear controllers that render the mapping $r\to e$ finite gain $\ell^\infty$-stable.

\subsection{Optimization over all internally stabilizing controllers}\label{sec:LTIS}

Consider first the case where the optimization is performed over all internally stabilizing LTI controllers,  that is:
\[ C = \argmin_{C \text{int. stab.}} \left  \|\frac{1}{1+C(z)\frac{1}{z-a}}R(z) \right\|_2
\]
 Using the Youla parameterization \cite{SpSzBook} and expressing $S(z)$ in terms of the Youla parameter $Q(z)$ leads to the following (weighted) model matching problem:
\beq \label{eq:optimization}
 \min_{Q(z) \in \mathcal{RH}_\infty}\left \|R(z)M(z)(Y(z)+N(z)Q(z))\right\|_2 \eeq
 where $N$ and $M$ are a coprime factorization of the plant $P=\frac{1}{z-a}$ and where $Y(z)$ is a solution of the following Bezout equation in $X(z),Y(z) \in \mathcal{RH}_\infty$:
 \beq \label{eq:bezout1} N(z)X(z)+M(z)Y(z)=1 \eeq
 The problem above can be simplified by choosing a coprime factorization where $M(z)$ satisfies  $|M(z)|=1$ for all $|z|=1$  and thus, for all signals  $x \in \ell^2$, $\|M(z)X(z)\|_2=\|X(z)\|_2$:
 \begin{align}
    &N=\frac{1}{az-1}, \;  M=\frac{z-a}{az-1} \\
    &X=a^2-1, \; Y =a
\end{align}
leading to:
\beq \label{eq:reformulated}
\begin{aligned}
& \min_{Q \in \mathcal{RH}_\infty} \left \|R(z)M(z)[Y(z)+N(z)Q(z)] \right \|_2  = \\ 
&\min_{\tilde{Q} \in \mathcal{RH}_\infty} \left\|aR(z)+\frac{R(z)}{z}\tilde{Q}\right\|_2
 \end{aligned}
 \eeq
 where we used the fact that $M(z)$ is all pass and we defined 
\[\tilde{Q}(z) \doteq Q(z) \frac{z}{az-1} \] 
Problem \eqref{eq:reformulated} can be explicitly solved by considering an expansion of   the optimal $\tilde{Q}$ of the form 
\[
\tilde{Q}(z)= q_o + \frac{q_1}{z}+ \ldots \frac{q_{n-1}}{z^{n-1}} + \frac{\tilde{Q}(z)_{tail}}{z^n} 
\]
Parseval's Theorem, combined with the explicit expression for $R(z)$ yields:
\[\begin{aligned}
&\left  \|aR(z)+\frac{R(z)}{z}\tilde{Q} \right \|^2_2   =
r_o^2\left \| a +\frac{a+q_o}{z}+\frac{a+q_o+q_1}{z^2} \ldots+ \right. \\
& \left. \frac{a +q_o + \ldots q_{n-2}}{z^{n-1}} + \frac{\sum q_i}{z^{n}}+\frac{\mathcal{O}(q_1, \ldots,q_{n-1},\tilde{Q}_{tail})}{z^{n+1}}\right \|^2_2 \\ 
& = r_o^2\left [ a^2 +  (a+q_o)^2+ 
 \ldots+(a + q_o +\ldots q_{n-1})^2 \right ] \\ & + r_o^2 (\sum q_i)^2 + r_o^2\left \|\mathcal{O}(q_1, \ldots,q_{n-1},\tilde{Q}_{tail})\right \|^2_2 
\end{aligned}\] 
Hence the optimal solution is given by
$q_o=-a$, $q_i=0, \;i \geq 1$, $\tilde{Q}_{tail}=0$, with the corresponding $Q$, controller $C$, closed-loop sensitivity $S$, and optimal cost given by
\beq \label{eq:optimal} \begin{aligned}
Q=-\frac{a(az-1)}{z}, \quad &C=\frac{(a^2+a-1)z-a^2}{a(z-1)} \\
S=\frac{a(z-a)(z-1)}{(az-1)z}, \quad
&S(z)R(z)=\frac{ar_o(z-a)(z^n-1)}{(az-1)z^n} \\ \|S(z)R(z)\|_2&=ar_o\sqrt{2}
\end{aligned}
\eeq
Since the cost $ar_0\sqrt{2}$ in \eqref{eq:optimal} is optimal, any controller yielding a lower cost \textbf{cannot be internally stabilizing.} Note that
since $Q$ is stable and proper, the closed loop system must satisfy  $S(\infty)=1$ and $S(a)=0$. These interpolation conditions follow  from \eqref{eq:reformulated} and the fact that $M(a)N(a)=M(\infty)N(\infty)=0$. Indeed, problem \eqref{eq:optimization} can be recast as:
\beq
   \label{eq:with_int} \min_{S \in \mathcal{RH}_\infty}\|S(z)R(z)\|_2 \quad \text{s.t.} \; S(\infty)=1, \; S(a)=0 
\eeq

\subsection{Input/Output optimization over LTI controllers}\label{sec:io}

In this section we show that simply optimizing  $\|S(z)R(z)\|_2$ without taking into account the interpolation constraint
$S(a)=0$ leads to controllers that are not internally stabilizing. Consider a controller in the form:
\beq \label{eq:controller_poles_zeros}
C_1(z)=K\frac{\prod_{i=1}^{n_z}(z-z_i)}{\prod_{i=1}^{n_p}(z-p_i)}. 
\eeq
Direct minimization of $\|S(z)R(z)\|_2$ with respect to $K,z_i,p_i$, with $n_p=n_z=1$, using Matlab's \cite{MATLAB} command {\tt fminsearch}\footnote{Interestingly, this model agnostic optimization yields a controller with a pole at $z=1$, which is consistent with the internal model principle \cite{IMP}.} leads to
\beq \label{eq:controller}
C_1(z) = \frac{(z-a)}{(z-1)}\eeq
that yields the closed--loop mapping $r\to e$:
\beq \begin{aligned}
    S&=\frac{z-1}{z} \;\; \text{with} \;\;
    \|S(z)R(z)\|_2= \left\|\frac{z-1}{z}\cdot\frac {r_o(z^n-1)}{(z-1)z^{n-1}} \right\|_2 \\ &= r_o\|1 -\frac{1}{z^n}\|_2  = r_o\sqrt{2}
< ar_o\sqrt{2}\end{aligned}\eeq
By parameterizing all controllers as 
$C=Q(z-a)/z(1-\frac{Q}{z})$ and dropping the no unstable pole/zero cancellation requirement, it can be shown that this controller with $Q=1$ is globally optimal over the set of all LTI controllers that optimize $\|S(z)R(z)\|_2$ s.t. the input/output stability constraint.

\begin{figure}[hbt]
    \centering
    \begin{minipage}{0.49\linewidth}
        \centering
        \includegraphics[width=\linewidth]{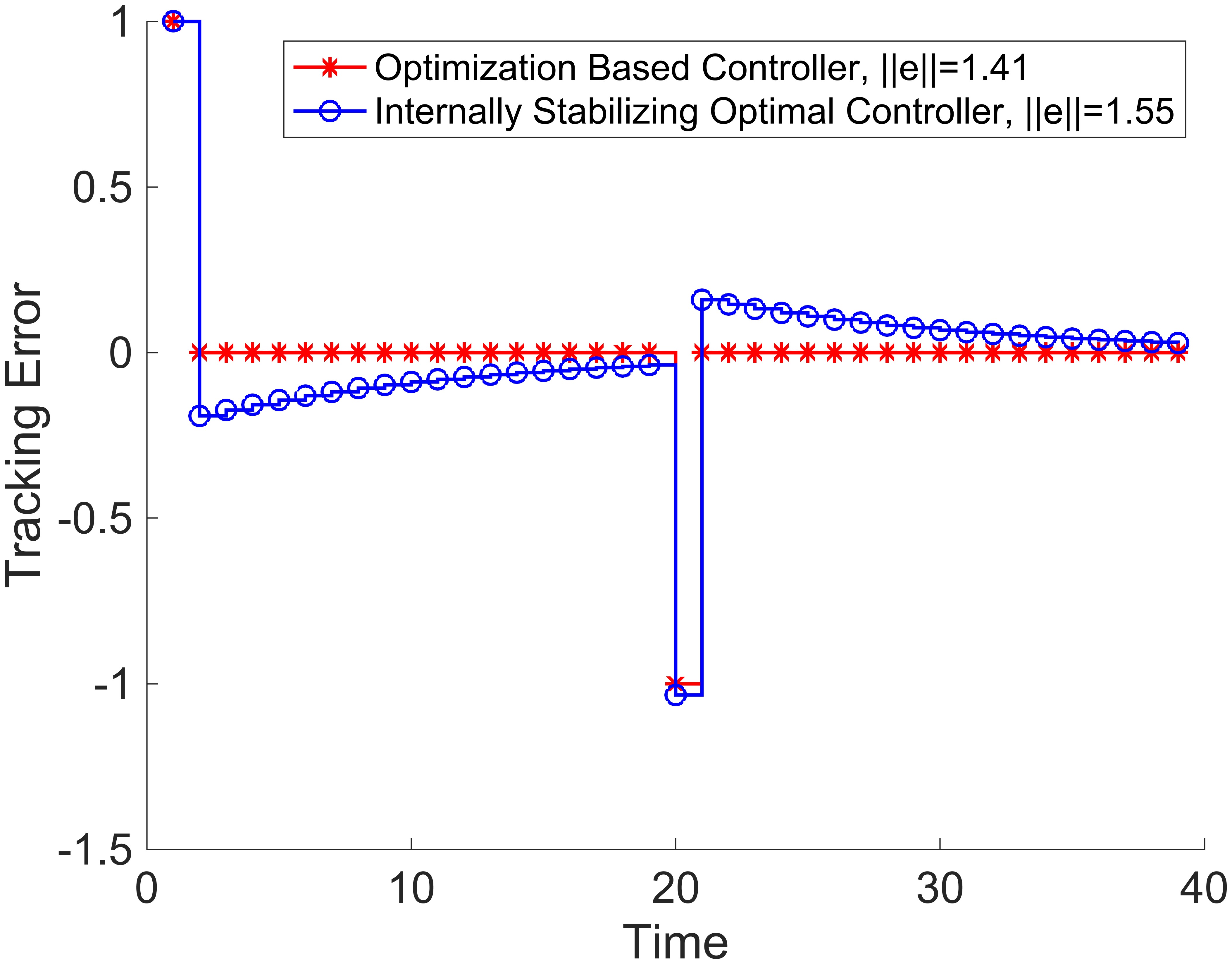}
        \\ (a)
    \end{minipage}
    \begin{minipage}{0.49\linewidth}
        \centering
        \includegraphics[width=\linewidth]{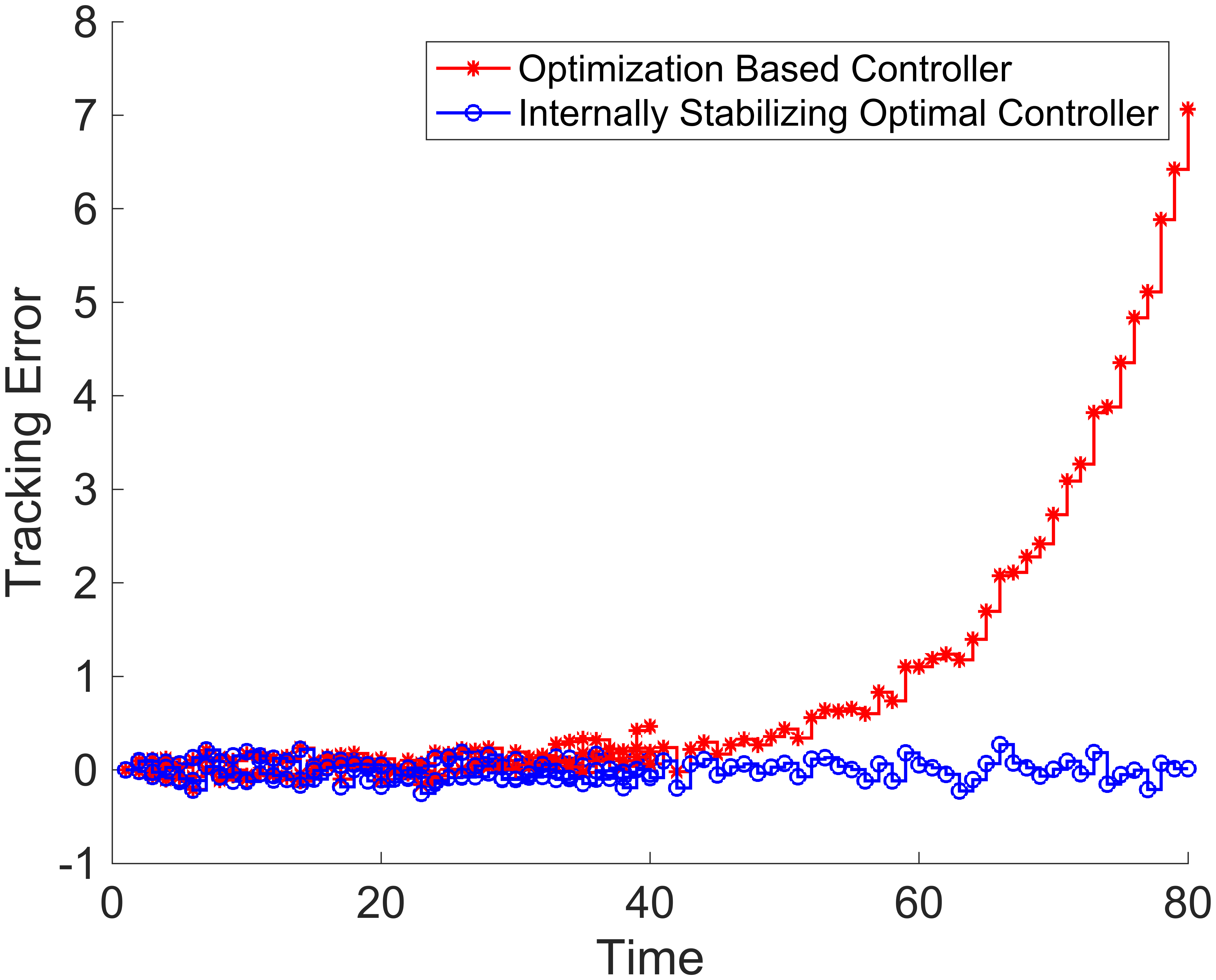}
        \\ (b)
    \end{minipage}
    \vskip -0.5 em
    \caption{Closed loop responses  for the controllers \eqref{eq:optimal} and \eqref{eq:controller}: (a) Tracking error; (b) Response to a random perturbation $w$. } \label{fig:comparison} 
\end{figure}
\vskip 0.7 em

A comparison of the closed loop response obtained using the controllers \eqref{eq:optimal} and \eqref{eq:controller} is shown in Fig. \ref{fig:comparison}. For the experiments, we use $r_o=1$, $a=1.1$. As expected, the controller \eqref{eq:controller} achieves a
lower  cost than \eqref{eq:optimal}.  However, the resulting closed loop system $r\to e$ is not internally stable, due to the unstable pole/zero cancellation. 
Specifically, the closed-loop transfer functions from $w$ to $e$ and $u$ are: 
\beq
T_{ew}=\frac{z-1}{z(z-a)}, \; T_{uw}=\frac{1}{z}
\eeq
Therefore, any arbitrarily small perturbation to the control action will lead to an unbounded output (Fig \ref{fig:comparison}, (b)). This instability does not show in the performance index being optimized. Thus, any algorithm that seeks to optimize it with respect to  the parameters of a controller of the form \eqref{eq:controller_poles_zeros} with $n_z \geq 1, n_p \geq 1$ and achieves global optimality will lead to a controller that has to perform as well as \eqref{eq:controller}, resulting in an input/output optimal controller that is not internally stabilizing.
Further, since $T_{uw}$ is stable, the control action remains bounded, even if the output does not. 
Hence the loss of internal stability cannot be detected by adding noise to the signal $r$ during training and monitoring the magnitude of the control.
Indeed, adding a disturbance $w \sim \mathcal{N}(0,0.1)$ to $r$ during training still leads to the controller \eqref{eq:controller} and the unstable pole/zero cancellation.

The flaw in the input/output optimization discussed above is that it does not enforce the interpolation conditions. Indeed, while the sensitivity $S=\frac{z-1}{z}$ satisfies the condition $S(\infty)=1$, it does not satisfy the interpolation condition $S(a)=0$. Removing the interpolation constraints leads to a super-optimal controller that is not internally stabilizing due to the unstable pole/zero cancellation. 

\begin{remark} The loss of internal stability  cannot be avoided by regularizing the performance index by adding  a penalty in the control action.  This penalty will avoid using integral action but still leads to an unstable pole/zero cancellation. For instance, changing the objective in the optimization above to $J=\|S(z)R(z)\|_2+\|u\|_2$ leads, for the case $a=1.1$, to
\[C(z)=\frac{0.60843 (z-1.1)}{
    z-0.9935}\]
which still exhibits the unstable pole/zero cancellation.
\end{remark}

\subsection{General Optimization based nonlinear controllers}

\begin{figure}[H]
    \centering
    \vskip -0.5 cm
    \includegraphics[width=0.9\linewidth]{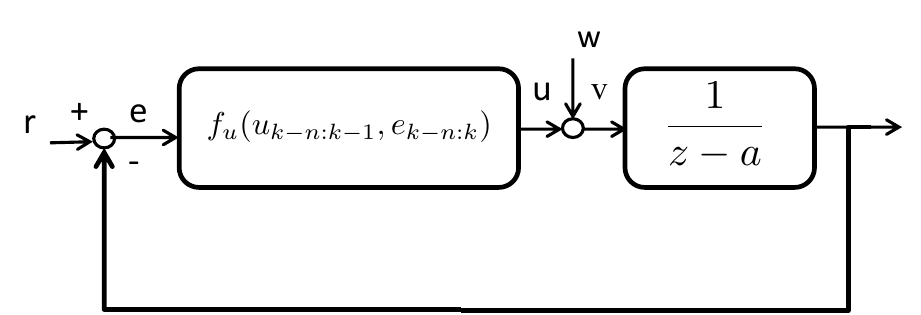} 
    \caption{A   controller that minimizes $\frac{\|e\|_2}{|r_o|}$  cannot render the mapping $\left[r \quad w\right]^T \to \left[u \quad e\right]^T$ finite $\ell^\infty$ gain  stable.}  \label{fig:CL3} 
\end{figure}
\vskip 0.5 em

In this section we show that any continuous nonlinear controller that (i) renders the input/output mapping $r \to s \doteq \left[u \quad e\right]^T$ finite-gain $\ell^\infty$ stable,  and (ii) achieves  tracking performance $\|e_k\|_2 = r_o\sqrt{2}$ for all $r_o$, will not internally stabilize the loop. Specifically, we will show that the resulting  closed loop system shown in Fig. \ref{fig:CL3} cannot have a finite $\ell^\infty$ gain from $\left[r \quad w \right]^T \to \left[u \quad e \right]^T$.
\begin{theorem}\label{teo:mainlinf} Consider 
a finite dimensional nonlinear controller  of the form 
\beq
\label{eq:controller1} 
 u_{k}=f_u(\greekvect{\theta}_k)\eeq
where $\greekvect{\theta}_k \doteq \begin{bmatrix}
\mat{u}_{k-1} \\ \mat{e}_k\end{bmatrix}$,  
$\mat{u}_{k-1} \doteq \begin{bmatrix}u_{k-1}
\ldots u_{k-m}\end{bmatrix}^T$, $\mat{e}_{k} \doteq \begin{bmatrix}e_{k} \ldots e_{k-m}\end{bmatrix}$, $m$ is the memory of the controller, and  $f_u$ is continuous. Let $\Phi_{cl}$ denote the corresponding closed loop mapping from the input  $r$ to the output sequence $\left \{\begin{bmatrix}u_k & e_k\end{bmatrix}^T \right\}$.  If the controller \eqref{eq:controller1} is such that:
\begin{enumerate}
    \item[(i)] $\Phi_{cl}$ is finite $\ell^\infty$ gain stable, e.g. $\left \|\begin{bmatrix}u & e\end{bmatrix}^T \right \|_\infty \leq K_r \|r\|_\infty$ and 
    \item[(ii)] when $r$ is a width-n pulse of the form \eqref{eq:rectangular},
    $\|e\|_2 = r_o \sqrt{2}$.
    \end{enumerate}
      then  the closed loop mapping $\begin{bmatrix}r & w \end{bmatrix}^T \to \begin{bmatrix} u &e \end{bmatrix}^T$ does not have finite $\ell^\infty$ gain. 
\end{theorem}

\begin{proof}
Consider the controller \eqref{eq:controller1} and note that finite closed-loop $\ell^\infty$ gain of $\Phi_{cl}$ implies that $f_u(\mat{0})=0$.  Assume for now that $f_u$ is twice differentiable and  consider its  linearization around $\mathbf{0}$:
\beq \label{eq:linearization} 
\begin{aligned}
u_k = & \frac{\partial f_u  (\mat{0})}{\partial \greekvect{\theta}} \begin{bmatrix}\mat{u}_{k-1} \\ \mat{e}_k \end{bmatrix}  +
\frac{1}{2}\begin{bmatrix}\mat{u}^T_{k-1}  \mat{e}^T_k\end{bmatrix} \frac{\partial^2 f_u(\greekvect{\theta}_o)}{\partial {\theta_i} \partial \theta_j}\begin{bmatrix}\mat{u}_{k-1} \\ \mat{e}_k \end{bmatrix}
\end{aligned}
\eeq
By contradiction, assume  that the closed loop mapping $\begin{bmatrix}r & w \end{bmatrix}^T \to \begin{bmatrix} u &e \end{bmatrix}^T$  has finite $\ell^\infty$ gain $K_w$. We will show that under this assumption, the 
 linear controller
\beq \label{eq:linearA}
u_k = \mathcal{L}(\mat{u}_{k-1},\mat{e}_k) \doteq  \frac{\partial f_u  (\mat{0})}{\partial u} \mat{u}_{k-1} + \frac{\partial f_u (\mat{0})}{\partial e}  \mat{e}_k
\eeq
is internally stabilizing. 
 By construction, the control sequence generated by the linear controller \eqref{eq:linearA} is the same sequence generated by the nonlinear one in the presence of a fictitious disturbance
 \[\begin{aligned}\hat{w}_k & = -\frac{1}{2}\begin{bmatrix}\mat{u}^T_{k-1}  \mat{e}^T_k\end{bmatrix} \frac{\partial^2 f_u(\greekvect{\theta}_o)}{\partial {\theta_i} \partial \theta_j}
\begin{bmatrix}\mat{u}_{k-1} \\ \mat{e}_k \end{bmatrix} \\ 
&\text{
with } \;
|\hat{w}_k| \leq 0.5 \left\|
\frac{\partial^2 f_u(\greekvect{\theta}_o)}{\partial {\theta_i} \partial \theta_j}
\right\|_2 \left \|\begin{bmatrix}\mat{u}_{k-1} \\ \mat{e}_k \end{bmatrix}\right \|^2 
\end{aligned}\] 
Hence \beq \label{eq:wboundA} \begin{aligned} \|\hat{w}\|_\infty &   \leq  C_1 \left\|\begin{bmatrix}\mat{u} \\ \mat{e} \end{bmatrix}\right \|^2_\infty 
\leq C_1 K_w^2\left \|\begin{bmatrix} r \\ w \end{bmatrix}\right \|_\infty^2 
\end{aligned}\eeq
for some constant $C_1$ that depends on $m$, the memory of the controller and the norm of the Hessian. Since the plant is linear, the error generated by the linear controller $\mathcal{L}$ satisfies
 \beq \label{eq:linearerrorA}e_{\mathcal{L}}  = e_{\mathcal{NL}}+ e_{\hat{w}}\eeq
 where $e_{\hat{w}}$ denotes the error due effect of the fictitious perturbation, with
 \beq \label{eq:ewboundA}
 \|e_{\hat{w}}\|_\infty \leq K_w \left \| \begin{bmatrix} r \\ \hat{w} \end{bmatrix}\right \|_\infty
 \eeq
 Thus
 \[ \begin{aligned} & \left \|\begin{bmatrix}u_\mathcal{L} \\ e_\mathcal{L} \end{bmatrix} \right \|_\infty   \leq \left \|
 \begin{bmatrix}u_\mathcal{NL} \\ e_\mathcal{NL} \end{bmatrix}
 \right \|_\infty + \left\| \begin{bmatrix} \hat{w} \\ e_{\hat{w}}\end{bmatrix}\right\|_\infty \\
 &\leq K_w \left\| \begin{bmatrix} r \\ w\end{bmatrix}\right\|_\infty + \left\| \begin{bmatrix} \hat{w} \\ e_{\hat{w}}\end{bmatrix}\right\|_\infty 
 \end{aligned}\]
 This inequality, combined with \eqref{eq:wboundA} and \eqref{eq:ewboundA}, shows that if $\begin{bmatrix}r \\ w \end{bmatrix} \in \ell^\infty $, then  $\begin{bmatrix}u_\mathcal{L} \\ e_\mathcal{L} \end{bmatrix} \in \ell^\infty$. Since both the plant and the controller $\mathcal{L}$ are linear, this implies that the mapping
 $\begin{bmatrix}r \\ w \end{bmatrix} \to \begin{bmatrix}u_\mathcal{L} \\ e_\mathcal{L} \end{bmatrix}$  has all its poles in the open disk $|z|<1$ (e.g. it is bounded input bounded output (BIBO) stable).  To complete the proof, we need to show that the linear controller $\mathcal{L}$ achieves a tracking error $\|e\|_2 < ar_o\sqrt{2}$ to the input \eqref{eq:rectangular}.  Let $T_{ew}$ denote the closed loop mapping $w \to e $ achieved by the linear controller $\mathcal{L}$.
Since the nonlinear controller achieves an error $\|e_{\mathcal{NL}}\|_2 =r_o\sqrt{2}$ to the input \eqref{eq:rectangular}, from the open-loop optimization in section \ref{sec:io}, it follows that ${u_{\mathcal{NL}}}_k =0$ and $e_k=0\; \forall k\geq n+1$, which implies $\hat{w}_k = 0, \; k\geq n+2+m$. Hence $\|\hat{w}\|_2 \leq (n+2+m)\|\hat{w}\|_\infty$.
  From  \eqref{eq:wboundA}, \eqref{eq:linearerrorA} and assumptions (i) and (ii) we have:
\beq \label{eq:LperformanceA} \begin{split}
\|e_\mathcal{L}\|_2 \leq \|e_\mathcal{NL}\|_2 + 
\|T_{ew}\|_{\ell^2 \to \ell^2}\|\hat{w}\|_2 
\leq  \|e_\mathcal{NL}\|_2 + \\
\|T_{ew}\|_{\ell^2 \to \ell^2}C_2K_r^2\|r\|_2^2 
\leq r_o \sqrt{2} + \mathcal{O}(r_o^2) < a r_o \sqrt{2}
\end{split}
\eeq
 if $r_o$ is small enough. It follows that the LTI controller \eqref{eq:linearA}
internally stabilizes the loop and achieves $\|e\|_2 < ar_o \sqrt{2}$ which contradicts \eqref{eq:optimal}.

Consider now the general case where $f_u(.)$ is continuous but not necessarily smooth.  Since by assumption $f_u(.)$  renders the mapping $\begin{bmatrix}r \\w \end{bmatrix} \to \begin{bmatrix} u \\e \end{bmatrix}$ finite $\ell^\infty$ gain stable, then,  as long as $r,w$ are confined to a compact set, so are $u,e$. From Stone-Weierstrass theorem \cite{Cheney} it follows that $f_u$ can be uniformly approximated arbitrarily close in this set by a polynomial $p_u(u,e)$, with $p_u(0)=f_u(0)=0$.  Thus,  the effect of this approximation can be absorbed into $\hat{w}$, as another term with $|w_{ap}| \leq \epsilon$, in the proof that $\mathcal{L}$ is internally stabilizing. In terms of performance, since $p_u(0)=0$ by construction, a same reasoning as before shows that, for the input \eqref{eq:rectangular}, $w_{ap}=0$ for $k >n+m+2$. Hence, $\|w_{ap}\|_2 \leq (n+m+2)\epsilon$ and the proof that the controller $\mathcal{L}$ (obtained now by linearizing the polynomial $p_u$) achieves $\|e\|_2 < ar_o\sqrt{2}$ still holds.
\end{proof}

\begin{figure}[hbt]
    \centering
    \begin{minipage}{0.49\linewidth}
        \centering
        \includegraphics[width=\linewidth]{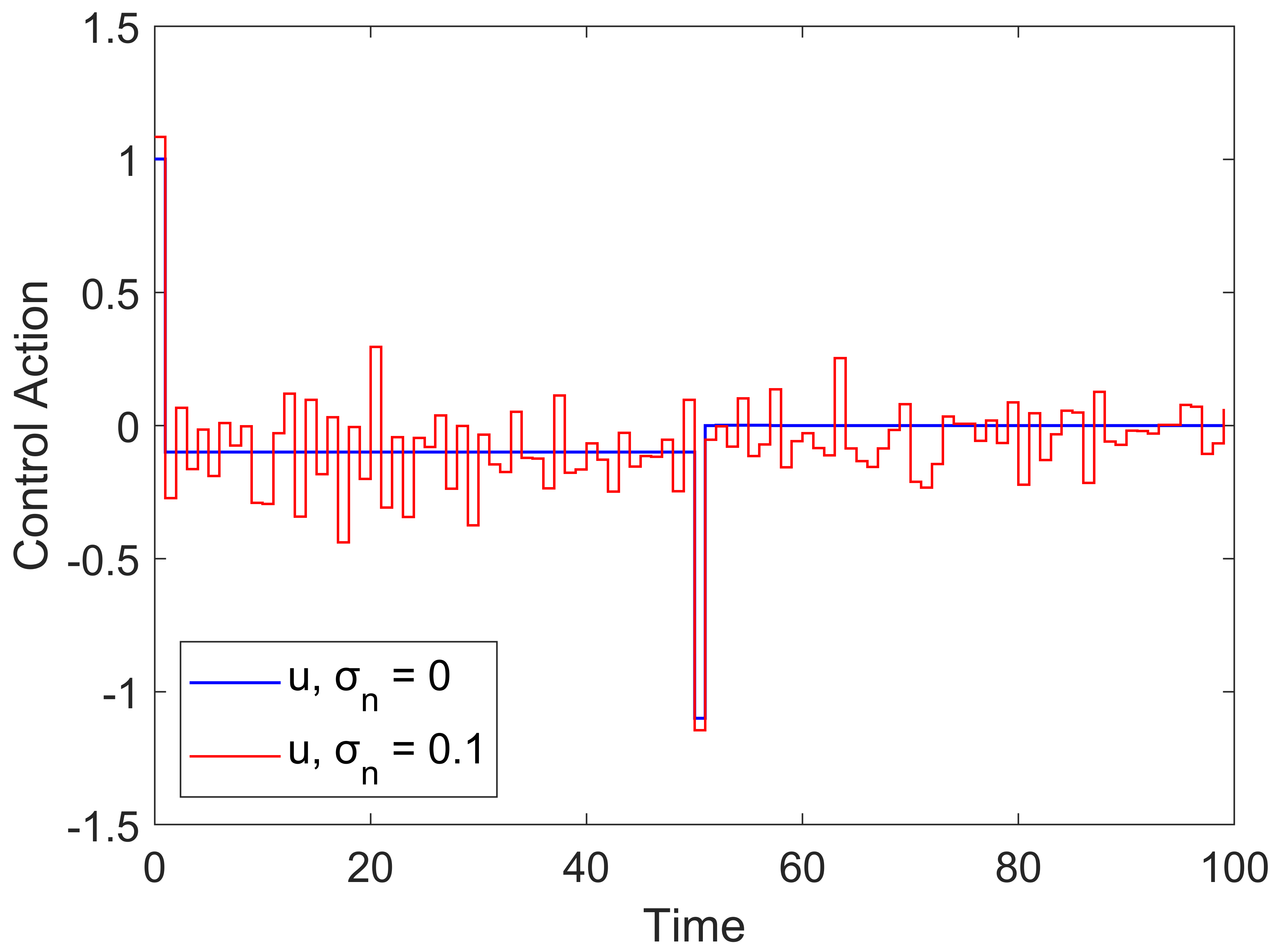}
        \\ (a)
    \end{minipage}
    \begin{minipage}{0.49\linewidth}
        \centering
        \includegraphics[width=\linewidth]{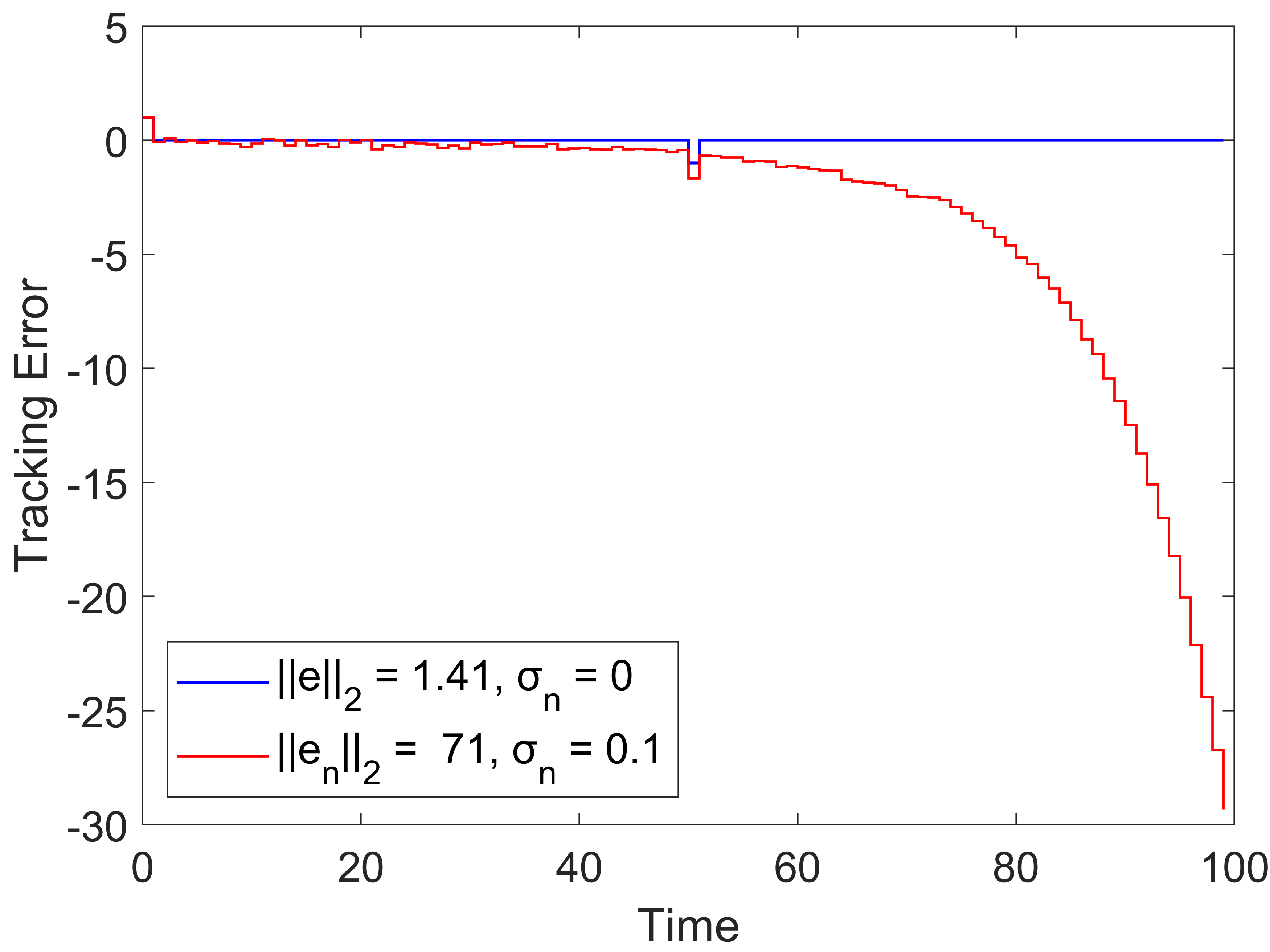}
        \\ (b)
    \end{minipage}
\vskip -0.5 em
\caption{Illustration of Theorem \ref{teo:mainlinf}. (a) Control action for the neural net controller, with and without noise added. (b) A small perturbation $w$ leads to an unbounded output.} 
\label{fig:nonlin_optim_e} 
\end{figure}

To empirically validate our results, we use a simple neural network controller with the following architecture:
\begin{equation}\label{eq:nn}
    \begin{aligned}
      &  \textbf{Input:} \; x \in \mathbb{R}^{3\times 1}, \; \textbf{Output:} \; y:=z^{(2)}, \quad y \in \mathbb{R}, \\
      &  \textbf{FC1:}\; z^{(1)} := W^{(1)}x + b^{(1)}, \; W^{(1)} \in \mathbb{R}^{2\times 3}, \, b^{(1)}\in \mathbb{R}^{2\times 1}\\
    &    \textbf{Activation:} \; a^{(1)} := ReLU(z^{(1)}), \; a^{(1)} \in \mathbb{R}^{2\times 1}\\
     &   \textbf{FC2:} \; z^{(2)} := W^{(2)}a^{(1)} + b^{(2)}, \; W^{(2)} \in \mathbb{R}^{1\times 2}, \, b^{(2)}\in \mathbb{R}
    \end{aligned}
\end{equation}
The input to the neural network is the vector $\big[ u[k-1], e[k-1], e[k] \big]$, which includes the previous values for the control action and tracking error to fit the described mapping in \eqref{eq:controller1}. The bias terms $b^{(1)}$ and $b^{(2)}$ are set to $b^{(1)}=[0,0,0]^T$ and $b^{(2)}=0$ due to the finite $\ell_2$ gain assumption. In this simple case, the optimal weights $W^{(1)}:=[-0.6460,0.7106,-0.6460;0.4119,-0.4335, 3.1555]$,
$W^{(2)}:=[-1.5480,0.3169]$ were found in 36 seconds using Matlab's command \texttt{fminsearch}  to minimize the  tracking error to a pulse of the form \eqref{eq:rectangular} with both positive and negative amplitudes and different pulse-widths. In Figure \ref{fig:nonlin_optim_e}, we observe the neural network achieving the optimal $\|e\|_2$ given in \eqref{eq:controller1} as $r_o \sqrt{2}$. However, as we introduce a zero-mean random normal noise signal with $\sigma_n=0.1$ to the control input $u$, as expected, the output becomes unbounded.

\section{Possible solutions to avoid losing internal stability and their limitations.}\label{sec:analysis}

As noted in Section \ref{sec:io}, the loss of internal stability can be traced to the existence of interpolation conditions. Direct unconstrained  optimization of  the tracking error leads to super-optimal controllers that violate these constraints. This issue can be solved by enforcing the interpolation constraints during the optimization. However, this can be difficult to  
accomplish in a model-agnostic setting. Below, we briefly discuss some options for enforcing these constraints.

\noindent \emph{Adding noise while training:} 
\begin{figure}[hbt]
\centering \vskip -1em
\includegraphics[width=0.75\linewidth]{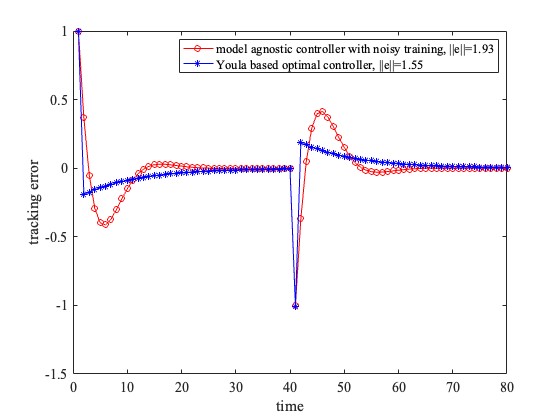} 
\vskip -1 em
\caption{Adding noise to the control action during training avoids the loss of internal stability at the price of 25\% performance degradation.} \label{fig:noisy_training} 
\end{figure}
Recall  that internal stability is equivalent  to input-output stability, provided that there are no unstable pole/zero cancellations between the plant and the controller \cite{SpSzBook}. Thus, perhaps the simplest way to implicitly avoid pole/zero cancellations is to add random noise to the control action during training, a technique that 
was empirically shown to improve robustness in \cite{Seiler2019}. 
While this is easy to implement in a model agnostic framework, it has the drawback of leading to suboptimal performance. This is illustrated in Fig. \ref{fig:noisy_training}, where adding a signal $w \sim \mathcal{N}(0,10^{-4})$ leads to an internally stabilizing controller, albeit with a 25\% performance degradation.  

\begin{figure}[hbt]
\centering
\includegraphics[width=0.7\linewidth]{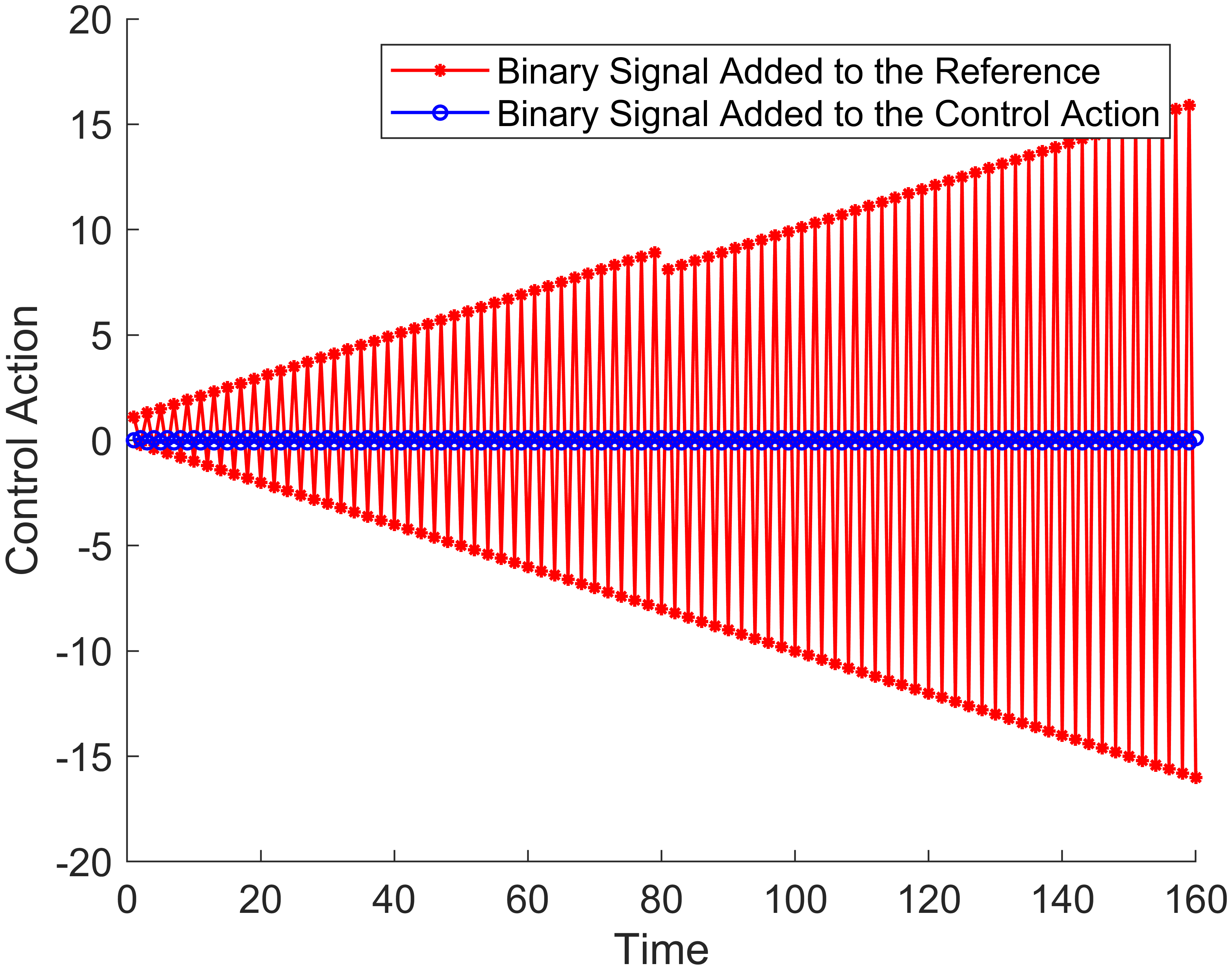} 
\vskip -0.5 em
\caption{Closed loop control response to perturbations in the reference (red) and the control (black), for a model agnostic controller for  the non-minimum phase plant $G=\frac{z+1}{z^2}$. Perturbing the reference signal leads to unbounded control.}
\label{fig:nonminphase} 
\end{figure}

A further problem is that  this approach critically hinges on adding noise at the right location. For 
instance, as illustrated in Section \ref{sec:io}, if the noise was instead added to the reference signal $r$, the resulting controller will still exhibit the unstable pole/zero cancellation.  Similarly, it can be shown that if the plant is changed to (the non-minimum phase one) $G=\frac{z+1}{z^2}$, then the optimal model agnostic controller is given by
$C=\frac{z^2}{(z+1)(z-1)}$. Fig \ref{fig:nonminphase} shows the control action for a model agnostic controller, trained with noise added to the control action, in response to perturbations $w_1$ and $w_2$, both chosen as a  binary signal with amplitude $0.1$, added to the reference input $r$ and the control $u$, respectively. As shown there, in this case the control action in response to $w_1$, the perturbation of the reference,  grows unbounded, even though the system was trained with noise added to the control. This highlights the importance of placing the perturbations at the ``right" place when training, so that the interpolation constraints are implicitly enforced.

\vskip 0.5 em

\begin{figure}[hbt]
\centering
\begin{minipage}{0.49\linewidth}
    \centering
    \includegraphics[width=\linewidth]{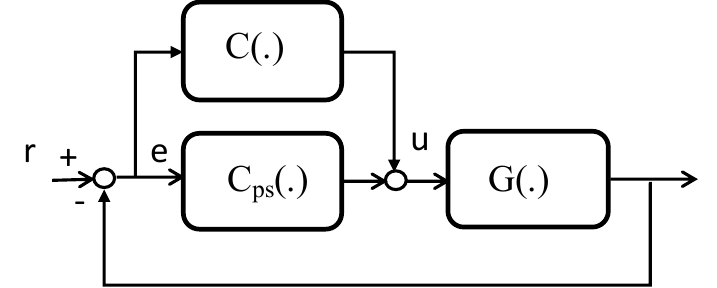}
    \\(a)
\end{minipage}
\hfill
\begin{minipage}{0.49\linewidth}
    \centering
    \includegraphics[width=1.1\linewidth]{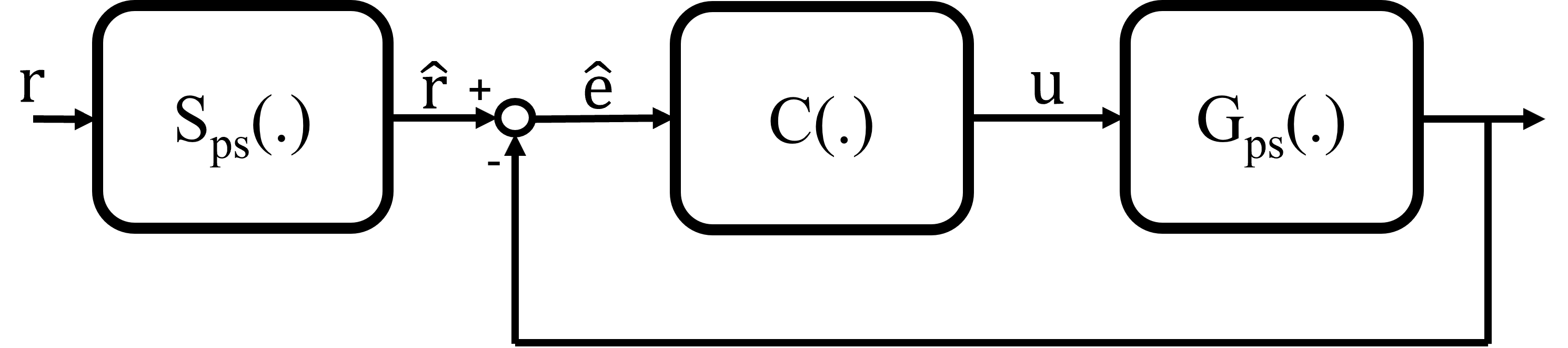}
    \\(b)
\end{minipage}
\caption{Prestabilizing a minimum-phase plant leads to optimal model agnostic controllers. (a) Loop showing the overall control action. (b) Equivalent  input/output mapping.}\label{fig:prestabilize} 
\end{figure}

\noindent \emph{Prestabilizing the plant:} An alternative to adding noise is to use a two step process where a prestabilizing controller $C_{ps}$ is learned first and then a second controller is added to optimize performance.
In this case, the mapping $r\to e$ is given by (see Fig \ref{fig:prestabilize})
\beq \label{eq:Ter} \begin{aligned}
T_{er}&=\frac{1}{1+G(C+C_{ps})}=S_{ps}\frac{1}{1+G_{ps}C}
\end{aligned}\eeq
 where
\beq \label{eq:prestab1} G_{ps} \doteq \frac{G}{1+GC_{ps}} \; \text{and} \; S_{ps} \doteq \frac{1}{1+GC_{ps}} \eeq

The advantage of this approach is that, as long as the unknown plant is minimum phase, the interpolation constraints (and hence internal stability) are automatically satisfied if the controller $C$ achieves input/output stability of the prestabililized plant, that is $\frac{1}{1+CG_{ps}}$ is stable. This follows from Lemma 2.3 in \cite{SpSzBook} and the fact that if $G$ is minimum phase so is $G_{ps}$. Since $G_{ps}$ is stable by construction, there cannot be unstable pole/zero cancellations between   $C$ and $G_{ps}$.

\vspace{4pt}

Further, in this case, the controller $C$ can recover any performance achieved by a controller working directly with the original plant. Since by assumption the plant $G$ is minimum phase, it satisfies the so-called parity-interlacing property. Hence, it can be stabilized with an open-loop stable controller (\cite{SpSzBook}, page 91). In turn, this implies that pre-stabilizing the plant does not affect achievable performance (\cite{SpSzBook}, page 87).

\begin{figure}[hbt]
    \centering
    \includegraphics[width=0.75\linewidth]{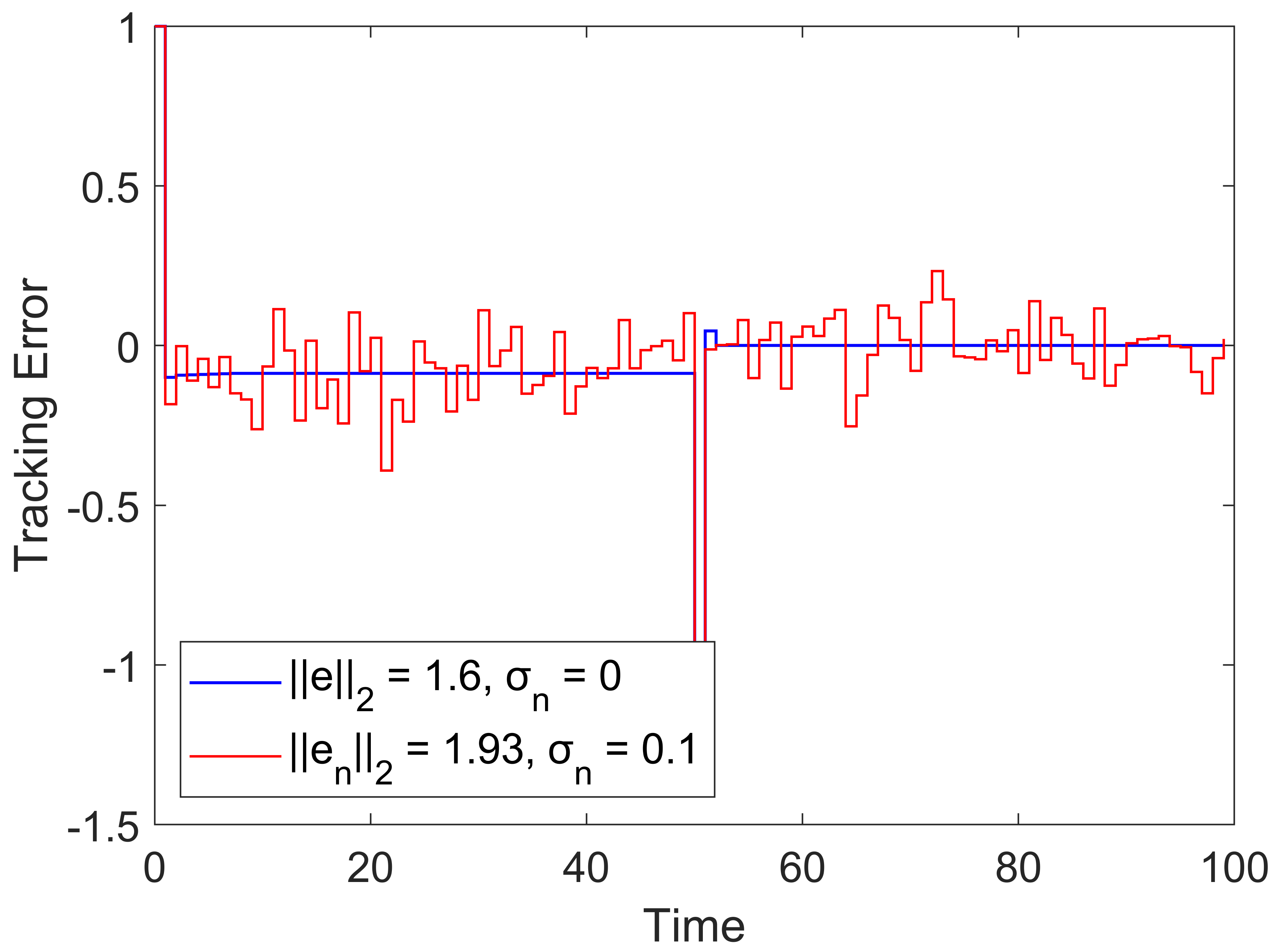}
    \vskip -0.5 em
    \caption{Tracking error $e$ achieved by the overall NN controller
    $C+C_{ps}$ with and without noise added to the control action.} \label{fig:gps_rhat} 
\end{figure}

For the simple example in this paper, one can search for a static stabilizing controller by simply minimizing the $\ell^2$ norm of the impulse response $\|1/(1+K/(z-a)\|_2$.
  Therefore, we use a 3 layer single-input, single output neural network with the input $e[k]$ to estimate such a controller. The resulting $C_{ps}$ has the weights $W^{(1)}:=[0.4561;1.5840;-0.7662]$, $W^{(2)}:=[0.1.2235,0.3421,-1.4357]$, with a ReLU activation in between and biases set to zero.  Once a suitable $C_{ps}$ has been found,  we then use the neural network in (\ref{eq:nn}) as the $C$ in Fig \ref{fig:prestabilize}(b) and optimize the weights. The resulting controller weights (computed in 428 seconds) are $W^{(1)}:=[0.5138,0.3661,-0.6312;2.3475,0.5361,-0.8748]$,
$W^{(2)}:=[-0.7884,0.5315]$. The tracking error $e$ achieved by the overall  controller ($C+C_{ps}$) is shown in Fig \ref{fig:gps_rhat}, both with and without noise added to the control action. As shown there, for the simple example in this paper, this approach indeed  leads to an internally stabilizing controller that achieves near optimal performance.

Drawbacks of this approach include the need for having a pre-stabilizing controller $C_{ps}$, which could be non-trivial to find, and its limitation to minimum phase plants.  In principle, non-minimum phase plants can be handled by adding noise to the control action, as outlined before, but this could entail performance loss. Further, if the plant is non-strongly stabilizable, the two-step approach may not be able to recover optimal performance.

\noindent \emph{Learning a coprime factorization of the  plant:}  While not strictly model agnostic, this is a data-driven approach where a model of the plant is learned first and then used to design a controller. Consider first the case of (unknown) LTI plants. The main idea  is  (i) to learn a coprime factorization of the plant, $G=NM^{-1}$, where the factors $N,M$ satisfy
$\eqref{eq:bezout1}$ for some stable $X,Y$, (ii) construct a prestabilizing controller $C_{ps}=X/Y$ (see \cite{SpSzBook}), and (iii) use the procedure outlined above.  Alternatively, one could directly use the parameterization of all stabilizing controllers 
$C=\frac{X-MQ}{Y+NQ}$ and optimize over the parameter $Q$. A potential difficulty here is that, due to the use of finite, noisy data records, only approximations $\tilde{N},\tilde{M}$ are learned. Thus, in principle there is no guarantee that the controller ${C}_{ps}$ will prestabilize the actual plant. As shown in a recent paper \cite{Rajiv2022}, this can be addressed by learning the factors $\tilde{N},\tilde{M}$ using as a loss function the gap metric between $(\tilde{N},\tilde{M})$ and $({N},{M})$. The advantages of using this metric (vis-a-vis other metrics) is that, if the resulting gap is below a quantity that can be directly computed from the experimental data, then the controller ${C}_{sp}$ is guaranteed to stabilize the true plant. Further, the factors $\tilde{N},\tilde{M}$ can be learned by solving a convex optimization problem.  The main disadvantage of this approach is that, at present time, computing the gap requires knowledge of the frequency response of the plant.  Hence, it cannot be directly applied to the case of interest in this paper where only finite time domain data is available. A time-domain characterization of the gap metric is needed to address this issue and to extend the approach to nonlinear plants.

\section{Conclusions}\label{sec:conclusions}
Model free direct policy optimization has the promise of optimizing performance while avoiding  a computationally expensive system identification step.  Further, it has been empirically shown to lead to optimal  controllers in a number of cases of practical importance. However, as we show in this paper with a very simple example, the success that direct policy optimization has achieved in some classes of problems is not directly generalizable to other seemingly simple problems, where it can lead to closed-loop systems that are fragile to arbitrarily small perturbations. This effect can be traced to the fact that completely model agnostic optimization can lead to unstable pole/zero cancellations and hence loss of internal stability. With this in mind, we proposed several ways to overcome this difficulty, provided that some minimal a priori information about the unknown plant is available.  Our results also point out to the need to develop a framework for learning stabilizing controllers  from finite, time domain data records, using as loss function the gap metric. This metric is, at its  core, a distance between closed-loop systems, as opposed to the  open-loop metrics more commonly used when learning models from data.

\appendix

\subsection{Control Theory background material}

To make the paper self-contained, next we recall basic results on Coprime Factorizations, Inner Transfer Functions, the Youla parameterization for SISO plants, and Interpolation Constraints. The interested reader is referred to Section 3.7 in the textbook \cite{SpSzBook} for details.

\subsubsection{Coprime Factorizations}\label{sec:coprime}

\begin{definition}
Two stable, proper discrete time transfer functions $N(z)$ and $M(z)$ are said to be coprime if there exist stable, proper transfer functions $X(z)$ and $Y(z)$ such that the following Bezout equation 
\beq \label{eq:Bezout}
N(z)X(z)+M(z)Y(z)=1
\eeq
\end{definition}
Intuitively, $N$ and $M$ are coprime if they do not have common zeros in $|z|\geq 1$.

Consider now a discrete time SISO plant with transfer function $G(z)$. It can be shown that $G(z)$ can always be written as
$G(z)=N(z)/M(z)$, where $M$ and $N$ are coprime. A specific realization for $N$ and $M$ for the SISO case is given below (see for instance in Lemma 3.7 in \cite{SpSzBook}).

Let 
\begin{eqnarray}\label{eq:plant_youla}
G(z) &\equiv &\left[\begin{tabular}{c|c}
     $A$&$B$ \\ \cline{1-2} $C$&$D$\end{tabular}\right]
\end{eqnarray}
be a stabilizable and detectable realization of $G(z)$ and $F$
and $L$ any matrices such that $A+BF$ and $A+LC$ are stable.
Define:
\beq
\label{rcfG}
\begin{bmatrix} M(z) & N(z) \end{bmatrix}
 \equiv \left[\begin{tabular}{c|c}
     $A+BF$&$B$ \\ \cline{1-2} $F$&$I$ \\
     $C+DF$&$D$\end{tabular}\right]
\eeq

Then (\ref{rcfG}) is a coprime
factorizations of $G(z)$ with the corresponding $X$ and $Y$ given by

\beq
\label{rcfX}
 \begin{bmatrix}X (z) & Y(z) \end{bmatrix}
 \equiv \left[\begin{tabular}{c|cc}
     $A+LC$&$-(B+LD)$&$L$ \\ \cline{1-3} $F$&$I$&$0$
     \end{tabular}\right]
\eeq

\subsubsection{Inner Transfer Functions}
\begin{definition} A stable, proper transfer function $T(z)$ is said to be inner if $T(z)T(\frac{1}{z})=1$.
\end{definition}    
The importance of inner transfer functions is that they preserve the $\ell^2$ norm, that is, for any stable $G(z)$ we have:
\[\|T(z)G(z)\|_2 = \|G(z)\|_2\]
This follows from the fact that
\[ \begin{aligned}
\|T(z)G(z)\|_2^2 = \frac{1}{2\pi j}\oint \frac{T(\frac{1}{z})G(\frac{1}{z})T(z)G(z)}{z}dz \\ =\frac{1}{2\pi j}\oint \frac{G(\frac{1}{z})G(z)}{z}dz = \|G(z)\|_2^2
\end{aligned}\]
\subsubsection{Youla Parameterization}\label{sec:Youla}
Consider a SISO plant $G(z)$ and its associated coprime factorization $G(z)=N(z)M^{-1}(z)$ Then all linear time-invariant controllers that internally stabilize the plant can be parametrized in terms of a free parameter $Q(z)\in \mathcal{RH}_\infty$ as:
\[C(z) = \frac{X-MQ}{Y+NQ}\]
In particular, setting $Q=0$ shows that the controller $C=X/Y$ is internally stabilizing.

The advantage of the Youla parameterization is that all relevant closed loop transfer functions become affine functions of $Q$ of the form 
\beq \label{eq:phi}
\Phi(z)=T_1(z)+T_2(z)Q(z)\eeq
where $T_1$ and $T_2$ depend on the open loop plant. For the specific transfer functions of interest in this paper we have
\begin{align}
  S(z)& =\frac{1}{1+G(z)C(z)} = M(z)(Y(z)+N(z)Q(z)) 
\end{align}
 Hence, optimizing a convex function $\mathcal{L(.)}$ of these transfer functions over the set of all stabilizing controllers becomes a convex optimization problem.  For instance, the problem of interest in this paper,  minimizing the weighted norm of $S(z)$  over all internally stabilizing controllers, reduces to a model matching problem of the form:
\[ \textrm{min}_{Q \in \mathcal{RH}_\infty} \|R(z)M(z)(1+N(z)Q(z))\|_2
\]
 In principle, the problem above is infinite dimensional, but in the case of interest to this paper, it admits a closed form solution. To see this, start by choosing the Youla parameterization so that $M(z)$ is inner (explicit formulas to achieve this are given in Section 6.6 in \cite{SpSzBook}). In that case the problem simplifies to
\[ \textrm{min}_{Q \in \mathcal{RH}_\infty} \|R(z)(1+N(z)Q(z))\|_2
\]
 As shown in Section \ref{sec:LTIS} using the explicit expressions for $R$ and $N$ leads to
\[ \begin{aligned}
Q&=-\frac{a(az-1)}{z}, \; C=\frac{(a^2+a-1)z-a^2}{a(z-1)} 
\end{aligned}\]
\subsubsection{Youla parameterization in the time domain}

\begin{figure}[hbt]
\centering

\includegraphics[width=0.5\linewidth]{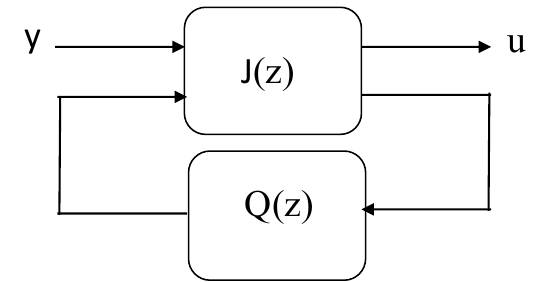} 

\caption{Parameterization of all stabilizing controllers as an LFT.}
\label{fig:LFT} 
\end{figure}

Consider the general output feedback problem for the plant
(\ref{eq:plant_youla}). Let $F$
and $L$ be constant matrices such that $A+LC_2$ and
$A+B_2F$ are stable. Then all the controllers that
internally stabilize $G$ are given by:
\beq
\label{allOF}
K(s)=F_{\ell}(J,Q), \qquad
 Q \in \mathcal{RH}_\infty, \; det \left [ I+D_{22}Q(\infty)\right ] \not = 0 
 \eeq
where
\beq
\label{JChap3}
J = \left[\begin{tabular}{c|c} $A+B_2F+LC_2+LD_{22}F$
&$\begin{array}{cc} -L &B_2+LD_{22}\end{array}$ \\
 \cline{1-2} $\begin{array}{c} F \\  -(C_2+D_{22}F)
 \end{array}$
 &$\begin{array}{cc} \mathbf{0} &I\\ I & -D_{22}
  \end{array}$\\ \end{tabular}\right]
\eeq
and where $F_{\ell}(J,Q)$ denotes lower Linear Fractional Transformation (see Fig. \ref{fig:LFT}).

\subsubsection{Parameterization of all achievable closed-loop maps.}

One of the main advantages of using the Youla parameterization is that   the resulting  closed--loop mapping is  affine in Q. This allows for recasting problems where the performance index is a convex function of input/output mapping  into a convex optimization in Q. Specifically, all achievable closed loop mappings have the form $F_{\ell}(T,Q)$
where $T =  \begin{bmatrix}T_{11} & T_{12} \cr T_{21} &
\mathbf{0}\end{bmatrix} $, or, equivalently 
\beq
T_{zw} =
T_{11}+T_{12}QT_{21}
\eeq
an affine function of $Q$. Here $T$ has the following state
space realization:
\beq
\label{AllTCh3}
\begin{array}{rcl}
T &=&\begin{bmatrix} T_{11} & T_{12} \cr
T_{21} & T_{22} \end{bmatrix}\\ \\
 &=&\left[\begin{tabular}{c|c} $\begin{matrix}A+B_2F &
-B_2F \cr \mathbf{0} & A+LC_2 \end{matrix}$
&$\begin{array}{cc} B_1 &B_2 \\ B_1+LD_{21} & \mathbf{0} \end{array}$ \\
 \cline{1-2} $\begin{array}{cc} C_1+D_{12}F  & -D_{12}F \\ \mathbf{0} & C_2 \end{array}$
 &$\begin{array}{cc} D_{11} & D_{12}\\ D_{21} & \mathbf{0}\end{array}$\\ \end{tabular}\right]
 \end{array}
 \eeq

In particular, setting $Q$=0 yields the familiar, separation based, controller structure corresponding to a state feedback gain $F$ and filter gain $L$.

\subsubsection{Internal Stability and the Interpolation Constraints}

Consider the generic optimal control problem of minimizing the norm of the closed loop system $\Phi(z)$ defined in \ref{eq:phi}. Let $z_i, i=1,\ldots n_z$ denote the non-minimum phase zeros of $T_2(.)$, that is those zeros of $T_2(.)$ with $|z_i|\geq 1$. Assume for simplicity that these zeros are simple. Clearly, since $Q(z)$ is stable (and hence has no poles in $|z|\geq 1$), a necessary condition for \eqref{eq:phi} to hold is that
\beq \label{eq:interp_constraints}
 \Phi(z_i)=T_1(z_i), \; i=1,\ldots, n_z
 \eeq
 Indeed, it can be shown that these conditions, known as the interpolation constraints, are necessary and sufficient for the existence of a stable $Q$ such that \eqref{eq:phi} holds (and hence internal stability). It follows that the problem of minimizing $\|\Phi\|_2$ over all internally stabilizing controllers can be recast into the following constrained optimization form:

\beq \label{eq:dual}
\begin{aligned}
   & \min_{\Phi \in \mathcal{RH}_\infty} \|\Phi(z)\|_2 \; \textrm{subject to \eqref{eq:interp_constraints}} 
\end{aligned}
\eeq

\subsection{Comparison with LQG control} 

 In principle the simple example shown in Fig. 2 can be recast into an $\mathcal{H}_2$ problem by  considering the generalized plant that has as inputs  an impulse disturbance and the control action, and as output the tracking error, that is
\[ \begin{bmatrix}e \\ e \end{bmatrix} =\begin{bmatrix} R(z)&  -\frac{1}{z-a} \\
 R(z)&  -\frac{1}{z-a}
 \end{bmatrix} \begin{bmatrix}\delta \\ u \end{bmatrix} \doteq G_p \begin{bmatrix}\delta \\ u \end{bmatrix} \]
 and minimizing the closed loop $\mathcal{H}_2$ cost from the input $\delta$ to the output $e$. However, the term $G_{21}$ has zeros on the unit circle at $z_k = e^{jk\frac{2\pi}{n}}$ and thus the standard two Riccati equations based algorithm fails (specifically, the Simplectic matrix associated with the filtering problem has eigenvalues on the unit circle, so a stabilizing solution to the ARE can't be found, please see section 2.1 of the textbook \cite{Kemin}).  In this case, the problem must be solved either using the Youla parameterization approach used in the paper, or by minimizing an upper bound of the $\mathcal{H}_2$ cost using LMI based arguments. 

\bibliographystyle{IEEEtran}
\bibliography{unstable_learning}

\end{document}